\documentclass{article}
 \newcommand {\theoremstyle} [1] { }

  \newenvironment{pot}{{\noindent\it\underline{Proof of Theorem \ref{nir}}}:}{\hfill$\Box$}
 \newenvironment{poc}{{\noindent\it\underline{Proof of Corollary \ref{grad}}}:}{\hfill$\Box$}
 \newtheorem{thm}{Theorem}[section]
 
 \theoremstyle{plain}
 \newtheorem{cor}{Corollary}[section]
 \theoremstyle{definition}

 \theoremstyle{remark}
 \newtheorem{rem}[thm]{Remark}

\usepackage {amssymb}
\usepackage {amsmath}
\usepackage{graphicx}

\def\R{\mathbb{R}}

\usepackage {amssymb}

\makeatother

\begin{document}

\author{Pablo Amster}

\title{{On Landesman-Lazer conditions and the Fundamental Theorem of Algebra} }
  \date{} 
  
\maketitle 

    \begin{abstract}

We give an elementary proof of a Landesman-Lazer type result for systems by means of a shooting argument and explore its connection with the fundamental theorem of algebra. 
    
    \end{abstract}

\section{Introduction}

In the well known paper \cite{landlaz}, Landesman and Lazer 
gave a sufficient condition for the existence of solutions of a nonlinear  
 scalar equation under resonance at a simple eigenvalue. 
 An extremely simplified version of this result is the periodic problem 
$$u'(t) + g(u(t))  = p(t),\qquad u(t+T)=u(t)
$$
where $g:\R\to\R$ is a smooth bounded function with limits $g_\pm$ at $\pm\infty$ and $p\in C(\R)$ is $T$-periodic. Here, the Landesman-Lazer condition  reads
$$g_-< \overline p < g_+$$
or
$$g_+ < \overline p < g_-$$
where $\overline p$ denotes the average of $p$, namely $\overline p := \frac 1T\int_0^T p(t)\, dt.$
Thus, the Landes-man-Lazer conditions express in fact two different things, that can be summarized as follows: 

 \begin{enumerate}
     \item $g_\pm \ne \overline p$. 
     \item The mapping $\Gamma:\{-1,1\}\to \R$ given by $\Gamma(\pm 1)=g_\pm$ wraps around $\overline p$, in the sense that 
          $(\Gamma(-1)-\overline p)$ and $(\Gamma(1)-\overline p)$ have different signs. 
     
 \end{enumerate}

In order to extend this idea for a system of differential equations, assume now that    $g:\R^2\to\R^2$  is bounded and the radial limits
$$g_v:= \lim_{r\to +\infty} g(rv),
$$
exist and are  uniform for $v\in \partial B$.
  Identifying  $\R^2$ with $\mathbb C$, we may define the curve $\Gamma(\theta):= g_{v(\theta)}$, with 
$v(\theta)=e^{i\theta}$ and $\theta\in [0,2\pi]$. The continuity of $\Gamma$ follows in a straightforward manner: for example, given $\varepsilon >0$ we may fix $r$ such that $|g(re^{i\theta}) - \Gamma(\theta)|<\frac \varepsilon 3$ for all $\theta$. Then
$$
|\Gamma(\theta) - \Gamma(\tilde\theta)| \le 
|\Gamma(\theta) - g(re^{i\theta})| + 
|g(re^{i\tilde\theta})-\Gamma(\tilde\theta)| + 
|g(re^{i\theta})-g(re^{i\tilde\theta})|.
$$
Using now the  continuity of $g$, there exists $\delta>0$ such that $|g(re^{i\theta})-g(re^{i\tilde\theta})|<\frac \varepsilon 3$, whence $|\Gamma(\theta) - \Gamma(\tilde\theta)|<\varepsilon$. 
In this setting, 
the following result due to Nirenberg \cite{nir} may be considered as a natural extension of the Landesman-Lazer theorem for a system. The condition that $\Gamma$ `wraps around' $\overline p$ shall be obviously expressed in terms of the winding number $I(\Gamma,\overline p)$:

\begin{thm}\label{nir}
In the previous situation, assume that
\begin{enumerate}
    \item $g_v\ne \overline p$ for all $v\in \partial B$. 
    \item $I(\Gamma,\overline p)\ne 0$. 
    \end{enumerate}
    Then the  problem 
\begin{equation}
\label{sys}    
u''(t) + g(u(t)) = p(t)
\end{equation}
has at least one $T$-periodic solution. 
\end{thm}

In the interesting paper  \cite{OS}, Ortega and S\'anchez observe that the Nirenberg condition does not hold 
for 
the so-called \textit{vanishing nonlinearities}, that is, when $g(u)\to \overline p$ as $|u|\to \infty$ and propose  to assume instead that $g(u) \ne \overline p$ for $|u|\gg 0$ and the limits
$$
q_v:=\lim_{r\to +\infty} \frac {g(rv) - \overline p}{|g(rv)-\overline p|}
$$
exist and are uniform for $v\in \partial B$. 
In this case,  Nirenberg's result is retrieved by defining now the (continuous) curve   $\Gamma_q(\theta):=q_{v(\theta)}$.

\begin{thm}\label{ort}

In the previous context, assume that
\begin{enumerate}
    \item \label{ort1}$q_v\ne \overline p$ for all $v\in \partial B$. 
 \item \label{ort2} $I(\Gamma_q,0)\ne 0$. 
    \end{enumerate}
    Then the system (\ref{sys}) has at least one $T$-periodic solution. 
\end{thm}

\begin{rem}
Observe that the second condition is analogous to the second condition in Theorem \ref{nir} due to the obvious fact that $I(\Gamma,\overline p)=I(\Gamma - \overline p,0)$.

\end{rem} 

  As a corollary,  it follows that if (\ref{sys}) is a  \textit{gradient system}, that is  
  $$
u'(t)=\nabla G(u(t)) + p(t),
$$ 
then the condition that $g=\nabla G$ is bounded can be dropped. The reason of this, as we shall see, is the fact that if $u$ is a $T$-periodic solution then multiplying the system by $u'(t)$ it is obtained, upon integration:
$$\int_0^T |u'(t)|^2\, dt = \int_0^T (G\circ u)'(t)\, dt + \int_0^T \langle p(t),u'(t)\rangle\, dt
$$
whence
$$\|u'\|_{L^2} \le \|p\|_{L^2}.
$$
\begin{cor} 
\label{grad}
Assume that $g=\nabla G$ and that conditions \ref{ort1} and \ref{ort2} of Theorem \ref{ort} are satisfied.  
    Then the system (\ref{sys}) has at least one $T$-periodic solution. 
\end{cor}

A particular instance of Corollary \ref{grad}
is the complex equation   
\begin{equation}\label{compl}
z'(t) = f(\overline{z}(t)) + p(t),    
\end{equation}
where $f$ is a polynomial. Indeed, in this case the Ortega-S\'anchez condition follows trivially since
$$\lim_{r\to+\infty} \frac{f(re^{-i\theta}) - \overline p }{|f(re^{-i\theta}) - \overline p|} = \frac{a_n}{|a_n|} e^{-in\theta},
$$
uniformly on $\theta$, where $a_n$ is the leading coefficient of $f$. This implies that $\Gamma_q$ performs $n$ clockwise turns around the origin and the conditions \ref{ort1} and \ref{ort2} in Theorem \ref{ort} are fulfilled. 

The fact that (\ref{compl}) is a gradient system follows from the  Cauchy-Riemann conditions: 
if $f=a+ib$ and $F= A + iB$ is a (complex) primitive of $f$, then 
$$[A(\overline z)]_x = A_x(\overline z) = a(\overline z), \qquad 
[A(\overline z)]_y = -A_y(\overline z) = b(\overline z).
$$
Alternatively, we may multiply the equation by $\overline z'(t)$ to obtain
$$|z'(t)|^2= z'(t)\overline z'(t) = G(\overline{z}(t))' + p(t)\overline z'(t);
$$
thus, if $z$ is a $T$-periodic solution it follows, as before, that
\begin{equation}\label{cotas}
\|z'\|_{L^2}\le \|p\|_{L^2}.
\end{equation}
This explains why 
a general version of the preceding result is interpreted by Mawhin 
in \cite{maw} as an extension of the Fundamental Theorem of Algebra: indeed, taking  
$p=0$, the inequality   (\ref{cotas}) implies that the periodic solutions are constants and, consequently, roots of $f$.

\section{Proofs}

In order to give elementary proofs of the preceding results, it proves convenient to recall a useful property of the winding number, which follows straightforwardly from the homotopy invariance: if $F:\overline {B_r(0)}\to\R^2$ is continuous and $I(\gamma,0)\ne 0$, where $\gamma(\theta):=F(re^{i\theta})$, then $F$ vanishes in $B_r(0)$. 

\medskip

\begin{pot}
Without loss of generality, we may assume that $\overline p=0$. 
Let $u(t)$ be a solution of (\ref{sys}) with initial value $u(0)=u_0$, then 
$$|u(t)-u_0| = \left| \int_0^t (g(u(s))+p(s))\, ds\right| \le M:= T(\|g\|_\infty + \|p\|_\infty).
$$
This implies that the Poincar\'e map $u_0\mapsto P(u_0):=u(T)$ is well defined, continuous and
$$P(u_0)-u_0= \int_0^T g(u(t))\, dt.
$$
Writing  $u_0= re^{i\theta}$ with $r>0$, it follows that $u(t)= r[e^{i\theta} + a(t)]$ with $|a(t)|\le \frac Mr$ and hence
$$u(t) = r(t) e^{i\theta(t)} 
$$
 where
 $$|r(t)-r|\le M, \qquad |\theta(t)-\theta|\le \frac Mr.
 $$
Thus, given $\varepsilon>0$, for sufficiently large $r$ we obtain
$$|g(u(t)) - \Gamma(\theta)|\le 
|g(u(t)) - \Gamma(\theta(t))| + |\Gamma(\theta(t))-\Gamma(\theta)| < \varepsilon
$$
and hence
$$|P(u_0)-u_0 - T\Gamma(\theta)| \le 
\int_0^T |g(u(t)) - \Gamma(\theta)|\, dt < T\varepsilon.
$$
Choosing   $\varepsilon < |\Gamma(\theta)|$ for all $\theta$ and setting $\gamma(\theta):= P(re^{i\theta}) - re^{i\theta}$, it follows that
$$
|\gamma(\theta) - \Gamma(\theta)| < |\Gamma(\theta)|
$$
for $r\gg 0$ which, in turn, implies that 
$$I(\gamma,0)= I(\Gamma,0) \ne 0.
$$
This proves the existence of $u_0$ such that $P(u_0)=u_0$, and the corresponding $u(t)$ is a $T$-periodic solution of the problem.  
\end{pot}
\medskip 

The proof of Theorem \ref{ort} is essentially the same as the preceding one: assuming w.l.o.g.  that $\overline p=0$, for $r\gg 0$ it is seen that 
$$I(\gamma_q,0) = I(\Gamma_q,0), 
$$
where $\gamma_q(\theta):= \frac{P(re^{i\theta}) - re^{i\theta}}{|g(r e^{i\theta})|}$, and the result follows. 

\medskip 
\begin{poc}
We may assume again that $\overline p=0$. 
With the aim of keeping the exposition at a very elementary level, let us assume for simplicity that 
$\nabla G$ is controlled by $G$, in the sense that 
\begin{equation}\label{control}
|\nabla G(u)|\le \xi(G(u))
\end{equation}
for some continuous mapping $\xi:\R\to (0,+\infty)$. For example, this holds   when $G$ is a polynomial, or if $G(u)=r(|u|)$ with $r\nearrow +\infty$. In this case, we may replace $G$ by a mapping $\hat G(u):=\varphi(G(u))$ with $\varphi:\R\to\R$ a smooth increasing function such that  
$$
\varphi'(s)=\left\{\begin{array}{cl}
1     &  \hbox{ if $|s|\le R$}\\
 \frac 1{\xi(s)}    & \hbox{ if $|s|\ge 2R$}
\end{array}
\right.
$$
for some $R$ to be specified. Observe that $\nabla \hat G(u) = \varphi'(G(u))\nabla G(u)$ is bounded and $\frac{\nabla \hat G(u)}{|\nabla \hat G(u)|} = \frac{\nabla G(u)}{|\nabla G(u)|}$, so by Theorem \ref{ort} the problem $u'(t)=\nabla \hat G(u(t))+p(t)$ has a $T$-periodic solution $u$. 
We claim that if $R$ is large enough, then $\|u\|_\infty\le R$ and, consequently, $u$ is a solution of the original problem. 
Indeed, as in the introduction it is verified that
 $$\|u'\|_{L^2} \le \|p\|_{L^2}
$$
and hence
$$|u(t)-u(0)| = \left| \int_0^t u'(s)\, ds\right|
\le T^{1/2}\|p\|_{L^2}:=M.$$ 
As before, fix $\varepsilon < |\Gamma_q(\theta)|$ for all $\theta$ and $r_0$ such that if $r\ge r_0$ then $$\left|\frac {\nabla G(r(e^{i\theta}+a))}{|\nabla G(re^{i\theta})|}-\Gamma_q(\theta)\right|<\varepsilon$$
for $|a|\le \frac Mr$. Because 
$\int_0^T \nabla \hat G(u(t))\, dt =0$, if $u_0=re^{i\theta}$ with
$r\ge r_0$ then we deduce that 
$$
\Gamma_q(\theta)\int_0^T \varphi'(G(u(t))\, dt= \int_0^T \varphi'(G(u(t)) \left[\Gamma_q(\theta)-\frac {\nabla G(u(t))}{|\nabla G(re^{i\theta})|}\right]\, dt.
$$
Thus
$$
|\Gamma_q(\theta)|\int_0^T \varphi'(G(u(t))\, dt < \varepsilon 
\int_0^T \varphi'(G(u(t))\, dt,
$$
a contradiction. Notice that $r_0$ depends only on $G$ and $M$; thus, it suffices to take $R=r_0+M$. 

\end{poc}

 \section{Further comments}
 
 It is easy to see that Theorem \ref{nir}  still holds   when $p$ is a  bounded function depending also on $u$; however, one needs to guarantee that, for $r$ large, the curve $\gamma(\theta):=P(re^{i\theta}) - re^{i\theta}$ wraps around $\frac 1T\int_0^T p(t,u(t))\, dt$, which varies  also with $\theta$. This is achieved if for example we assume
 \begin{equation}\label{limsup}
 \limsup_{|u|\to \infty} \frac{|p(t,u)|}{|g(u)|} = 0
     \end{equation}
  uniformly on $t$. It is readily verified that the same condition suffices also in the situations of Theorem \ref{ort} and Corollary \ref{grad}, assuming now that the limits
   \begin{equation}\label{radial}
 q^0_v:= \lim_{r\to+\infty} \frac {g(rv)}{|g(rv)|}       
   \end{equation} 
  exist uniformly for $v\in \partial B$ and replacing the curve $\Gamma_q$ by
  $$ \Gamma^0_q(\theta):= q^0_{v(\theta)} 
  $$
 The results may be extended also for delay systems like
 \begin{equation}
 \label{delay}
 u'(t)= g(u(t)) + p(t,u(t),u(t-\tau))
    \end{equation}
  where $\tau>0$ and $p$ is bounded, continuous and $T$-periodic in the first coordinate:   
  \begin{thm}
  Assume that $g$ is bounded or $g=\nabla G$ such that the radial limits (\ref{radial}) exist uniformly on $v\in \partial B$. Further, assume that $p$ is bounded with
     \begin{equation}\label{limsup-delay}
 \limsup_{|u|\to \infty} \frac{|p(t,u,u)|}{|g(u)|} = 0
     \end{equation}
     uniformly on $t$. 
If $I(\Gamma^0_q,0)\ne 0$, then  problem (\ref{delay}) has at least one $T$-periodic solution. 
  \end{thm}
It should be noticed that, in this case, the problem cannot be reduced to find a fixed point in a  finite dimensional space and less elementary tools are required. However, the proof is still easy in the context of the Leray-Schauder degree,  which yields the following continuation theorem:
\begin{thm}\label{LS}
Assume that
\begin{enumerate}
    \item There exists $R>0$ such that any $T$-periodic solution of the problem
    $$u'(t) =\lambda[g(u(t)) + p(t,u(t), u(t-\tau))]
$$
with $\lambda\in (0,1]$ satisfies $\|u\|_\infty <R$. 
\item $I(\Gamma_R,0)\ne 0$, where 
$$\Gamma_R(\theta):= g(re^{i\theta}) + \int_0^T p(t,re^{i\theta}, re^{i\theta})\, dt. 
$$
    
\end{enumerate}
Then  problem (\ref{delay}) has at least one $T$-periodic solution $u$ with $\|u\|_\infty< R$. 

\end{thm}

Indeed, when $g$ is bounded or $g=\nabla G$, it follows as before that 
  $\|u'\|_{L^2}$ is bounded by a constant depending only on $\|p\|_\infty$. This, in turn, implies that $|u(t)-u(0)|$ is also bounded and the conditions of Theorem \ref{LS} are obtained under the assumptions of Theorem \ref{ort}.  It is clear that  condition (\ref{control}) is not necessary at all. 
  
  Analogous  results may be obtained for larger systems: let $g\in C(\R^n,\R^n)$ and $p: \R\times \R^{2n}\to \R^n$ be continuous, bounded and $T$-periodic in its first coordinate. Assume that $g$ is bounded or $g=\nabla G$
   and that the radial limits (\ref{radial})
    exist uniformly for $v\in S^{n-1}\subset \R^n$. Furthermore, assume that (\ref{limsup-delay}) holds. Then the problem has at least one $T$-periodic solution, provided that the degree of the mapping $ \Gamma^0_{q}:S^{n-1}\to S^{n-1}$  given by $\Gamma^0_{q}(v):= q^0_{v}$ is different from $0$. 
   It is readily seen  that the latter condition is equivalent to 
   \begin{equation}\label{deg}
       \deg(g,B_R(0),0)\ne 0 \hbox{ when $R$ is sufficiently large,}
   \end{equation} 
   where $\deg$ stands for the Brouwer degree.  
  A more delicate argument given in \cite{ward} 
  shows that, if (\ref{deg}) is fulfilled, then the existence of the limits (\ref{radial}) is not necessary when 
   $g=\nabla G$ is coercive, that is $|\nabla G(u)|\to +\infty$ as $|u|\to+\infty$.

  
  As a final remark, let us try to understand why the result does not hold for the equation
  $$z'(t) = f(\overline z(t)) + p(t)
  $$
  when $f$ is an arbitrary entire function. Because $g$ is analytic, 
  the bounds for $z'$ are obtained exactly as before; however, if $f$ is not a polynomial then the curves $f(Re^{i\theta})$ with $R\gg 0$
are very badly behaved. This is clearly related  to the fact that $f$ has an essential singularity at $\infty$; for example, when  $f(z)=e^z$ the problem has no solutions for $p=0$ and  
 $$f(Re^{i\theta})= e^{Re^{-i\theta}} = e^{R\cos(\theta)} e^{-iR{\rm sin}(\theta)},
 $$
 which has zero winding number although for $R\gg 0$ it passes  many times back and forth around the origin.

\section*{Acknowledgements}
This research was  supported by the projects 
PIP 11220130100006CO CONICET
and UBACyT  20020160100002BA.

\bigskip 

\noindent 
Pablo Amster 

\noindent {Departamento de Matem\'atica. Facultad de Ciencias Exactas y Naturales\\
Universidad de Buenos Aires \& IMAS-CONICET\\
Ciudad Universitaria. Pabell\'on I,(1428), Buenos Aires, Argentina\\
\textit{{E-mail}}: pamster@dm.uba.ar}

  \end{document}